\documentclass[12pt]{article}
\usepackage{stmaryrd}
\usepackage{mathrsfs}
\usepackage{amsmath}
\usepackage{amsmath,amsthm,amssymb}
\usepackage{amsfonts}
\usepackage{latexsym}
\usepackage{CJK}
\date{}

\begin{document}
\title{Endomorphism algebras arising from mutations$^\star$}
\author{{\small Genhua Pei, Hongbo Yin, Shunhua Zhang}\\
{\small  School of Mathematics,\ Shandong University,\ Jinan 250100,
P. R. China }}

\pagenumbering{arabic}

\maketitle
\begin{center}
 \begin{minipage}{120mm}
   \small\rm
   {\bf  Abstract}\ \ Let $A$ be a finite dimensional algebra over an algebraically
   closed field $k$, $\mathcal {D}^b(A)$ be the bounded derived category of $A$-mod
   and $A^{(m)}$ be the $m$-replicated algebra of $A$. In this paper, we investigate
   the structure properties of endomorphism algebras arising from silting mutation in
   $\mathcal {D}^b(A)$ and tilting mutation in $A^{(m)}$-mod.
\end{minipage}
\end{center}

\vskip0.1in



\footnote{ $^\star$Supported by the NSF of China (Grant No.
11171183).}

\footnote{  Email addresses: \   again353535@163.com(G.Pei), \
yinhongbo0218@126.com(H.Yin),
 \ shzhang@sdu.edu.cn(S.Zhang).}

\vskip0.2in

\section {Introduction}

\vskip0.2in

Let $\mathscr{A}$ be an additive Krull-Schimidt category. It is well
known that the endomorphism algebras of rigid objects in
$\mathscr{A}$, in particular, of tilting modules over a finite
dimensional algebra have been central in representation theory
\cite{AR, BMR, BMRRT, HR, BZ}. In this paper, we focus on the
structure properties of endomorphism algebras arising from silting
mutation in $\mathcal {D}^b(A)$ and tilting mutation in
$A^{(m)}$-mod for a hereditary algebra $A$.

\vskip 0.2in

Let $A$ be a finite dimensional algebra over an algebraically closed
field $k$, and $\mathcal {D}^b(A)$ be the bounded derived category
of $A$-mod.  Let $T$ be a basic silting object in $\mathcal {D}$.
Without loss of generality we can assume that $T$ is in the
non-negative part $\mathcal {D}^+$  of $\mathcal {D}$, i.e. we
assume
$$
T = T_0[0] \oplus T_1[1]\oplus \cdots \oplus T_m[m] \ \ \ \ \ \ \ \
\ \ \ \ \ \ \ \ (1)$$ with each $T_i$ in ${\rm mod}~A$. We fix
$$
\mathcal {S}_m = {\rm mod}~A[0] \vee {\rm mod}~A[1]\vee \cdots \vee
{\rm mod}~A[m-1] \vee A[m].
$$
Then we prove the following theorems.

\vskip 0.2in

{\bf Theorem 1.}\ {\it Let $T$ be a basic silting object in
$\mathcal{S}_m$, and $\Gamma= {\rm End}_{\mathcal{D}}T$. Let
$\mathscr{T}= {\rm add}~T$ and $\mathcal {W}=\mathscr{T}\ast
\mathscr{T}[1]$. Then the functor
$$G = {\rm Hom}_{\mathcal {D}}(T, -
): \mathcal {W} \rightarrow {\rm mod}~\Gamma $$ is full and dense,
and $G$ also induces an equivalence functor
$$\overline{G} : \mathcal {W} /{\rm add}~(T[1]) \rightarrow  {\rm
mod}~\Gamma.$$}

\vskip0.2in

Let $T = T_1\oplus \cdots\oplus T_n$ be a basic silting object in
$\mathcal {S}_m$. Fix an indecomposable direct summand $T_i$, with
$1 \leq i \leq n$.  Recall from \cite{BRT}, we know that the almost
complete silting object $T/T_i$ has a countably infinite number of
non-isomorphic complements $M_i$ for $i\in \mathbb{Z}$. In
particular, there are complements $M_{-1}$ and $M_{m+1}$, such that
$M_j\simeq M_{-1}[j + 1]$ for $j < -1$, and $M_j = M_{m+1}[j-(m
+1)]$ for $j
>m +1$.

\vskip 0.2in

{\bf Theorem 2.}\ {\it Take the notations as above. Let $\Gamma_j=
{\rm End}_{\mathcal {D}}(M_j\oplus T/T_i)$ and $\mathcal
{W}_j=\mathscr{T}_j\ast \mathscr{T}_j[1]$, where $\mathscr{T}_j={\rm
add}~(M_j\oplus T/T_i)$. Then

\vskip 0.1in

(1)\ The  $\Gamma_j$ module ${\rm Hom}_{\mathcal {D}}(M_j\oplus
T/T_i, M_{j-1}[1])$ is simple.

\vskip 0.1in

(2)\  Let  $S_{M_j}$ be the simple top of the indecomposable
projective $\Gamma_j$ module ${\rm Hom}_{\mathcal {D}}(M_j\oplus
T/T_i, M_{j})$. Then we have
$$
{\rm mod}~\Gamma_j /{\rm add}~S_{M_j} \simeq \mathcal {W}_j/{\rm
add}~( M_{j}\oplus M_{j-1}\oplus T/T_i)[1].
$$
Moreover, if $j>m$, then $\Gamma_j \simeq \left (
\begin{array}{cc}
{\rm End}_{\mathcal {D}}M_j & 0\\
0 & {\rm End}_{\mathcal {D}}(T/T_i) \end{array}\right )$.}

\vskip 0.3in

Let $A^{(m)}$ be the $m$-replicated algebra of $A$. Recall from
\cite{Z},  we know that the tilting quiver of $A^{(m)}$ is
connected. For their endomorphism algebras, we have following
theorems.

\vskip 0.2in

{\bf Theorem 3.} \ {\it Let $M$ be a faithful almost complete
tilting $A^{(m)}$-module with non-isomorphic indecomposable
complements $X_0, \cdots, X_t$, and let $\Gamma_i={\rm
End}_{A^{(m)}}(X_i\oplus M)$. Then there is a BB-tilting $\Gamma_i$
module $T_i$ such that ${\rm End}_{\Gamma_i} T_i\simeq \Gamma_{i+1}$
for $0\leq i \leq t-1$.}

\vskip 0.2in

{\bf Theorem 4.}\ {\it Let $T=T_1\oplus T_2\oplus\cdots\oplus T_n$
be an basic silting object in $\mathcal{S}_m$, and let  $P$ be the
direct sum of all indecomposable injective-projective $A^{(m)}$
modules. Then $\Gamma_T= {\rm End}_{\mathcal{D}}~T\simeq {\rm
\overline{End}}_{A^{(m)}}~(T\oplus P)$.}

\vskip 0.2in

This paper is arranged as follows. In Section 2, we collect
definitions and basic facts needed for our research. Section 3 is
devoted to the proof of Theorem 1 and Theorem 2, and in section 4,
we prove Theorem 3 and Theorem 4.

\vskip 0.2in

\section {Preliminaries}

\vskip 0.2in

Let $\Lambda$ be a finite dimensional algebra over an algebraically
closed field $k$.  We denote by $\Lambda$-mod the category of all
finitely generated left $\Lambda$ modules, and by
$\Lambda$-$\overline{{\rm mod}}$ (resp. $\Lambda$-$\underline{{\rm
mod}}$) the factor category $\Lambda$-mod/$\mathscr{I}$ (resp.
$\Lambda$-mod/$\mathscr{P}$). The derived category of bounded
complexes of $\Lambda$-mod is denoted by $\mathcal{D}^{b}( \Lambda)$
and the shift functor by [1]. The positive part of $\mathcal{D}^{b}(
\Lambda)$ is denoted by $\mathcal{D}^{+}( \Lambda)$.

\vskip 0.2in

For a $\Lambda$ module $M$, we denote by ${\rm add}~M$ the
subcategory of $\Lambda$-mod whose objects are the direct summands
of finite direct sums of copies of $M$ and by
$\Omega_{\Lambda}^{-i}M$ the $i^{{\rm th}}$ cosyzygy
 of $M$.  The projective dimension of $M$ is denoted by
pd $M$, the global dimension of $\Lambda$ by gl.dim $\Lambda$ and
the Auslander-Reiten translation of $\Lambda$ by $\tau_\Lambda$.

\vskip 0.2in

Let $T$ be a $\Lambda$ module.  $T$ is said to be rigid if ${\rm
Ext}^i_\Lambda(T, T)=0$ for all $i\geq 1$.  A rigid module $T$ is
called a partial tilting module provided ${\rm pd} \ T<\infty$. A
partial tilting module $T$ is called a tilting module if  there
exists an exact sequence
$$0\longrightarrow \Lambda\longrightarrow T_0\longrightarrow T_1
\longrightarrow \cdots\longrightarrow T_d\longrightarrow 0$$ with
each $T_i\in {\rm add}\  T$. A partial tilting module $T$ is called
an almost complete tilting module if there exists an indecomposable
$\Lambda$-module $N$ such that $T\oplus N$ is a tilting module.

\vskip 0.2in

Let $\mathcal {A}$ be an additive category, and let $\mathcal{C}$ be
a full subcategory of $\mathcal {A}$, $C_{M}\in\mathcal{C}$ and
$\varphi :C_M\longrightarrow M$ with $M\in\mathcal {A}$. The
morphism $\varphi$ is a right $\mathcal{C}$-approximation of $M$  if
the induced  morphism ${\rm Hom}(C,C_{M})\longrightarrow {\rm
Hom}(C,M)$ is surjective for any $C\in\mathcal{C}$. A minimal right
$\mathcal{C}$-approximation of $M$ is a right
$\mathcal{C}$-approximation which is also a right minimal morphism,
i.e., its restriction to any nonzero summand is nonzero. The
subcategory $\mathcal{C}$ is called contravariantly finite if any
module $M\in\mathcal {A}$ admits a (minimal) right
$\mathcal{C}$-approximation. The notions of (minimal) left
$\mathcal{C}$-approximation and covariantly finite subcategory are
dually defined. It is well known that add $M$ is both a
contravariantly finite subcategory and a covariantly finite
subcategory. We call a morphism $\psi : X\longrightarrow Y$ in
$\mathcal{C}$  a sink map of $Y$ if $\psi$ is right minimal and
${\rm Hom} (\mathcal{C}, X)\longrightarrow {\rm Rad }(\mathcal{C},
Y)\longrightarrow 0$ is exact. A source map can be defined dually.

\vskip 0.2in

From now on, we denote by $A$ a finite dimensional hereditary
algebra over an algebraically closed field $k$.  The repetitive
algebra $\hat{A}$ of $A$ is the infinite matrix algebra
$$\hat{A} =\begin{pmatrix}
                                \ddots &  & 0 &  &  \\
                                 & A_{i-1} &  &  &  \\
                                & Q_{i} & A_{i} &  &  \\
                                &  & Q_{i+1} & A_{i+1} &  \\
                                & 0 & & \ddots \\
                             \end{pmatrix}
$$
which has only finitely many non-zero coefficients, $A_{i}=A$ and
$Q_{i}= DA$ for all $i\in \mathbb{Z}$, where $D={\rm Hom}_k(-,k)$ is
the dual functor, all the remaining coefficients are zero and
multiplication is induced from the canonical isomorphisms $A
\otimes_A DA\simeq \ _ADA_A\simeq DA\otimes_AA$ and the zero
morphism $DA\otimes_A DA \longrightarrow 0$ (see \cite{ABST, H,
HW}). Let
$$A^\backprime =\begin{pmatrix}
                  A_{0} &  &  &  & 0 & \\
                  Q_{1} & A_{1} & &  &  &  \\
                   & Q_{2} & A_{2} &  &  &  \\
                 &  &  & \ddots & \ddots &  \\
                                  \end{pmatrix}
$$
be the quotient of $\hat{A}$, which is called right repetitive
algebra of $A$.

\vskip 0.2in

{\bf  Lemma 2.1.}\ (1)$^{\cite{H}}$ \ {\it The derived category
$\mathcal{D}^{b}(A)$ is equivalent, as a triangulated category, to
the stable module category $\hat{A}$-${\rm \underline{mod}}$.}

\vskip 0.1in

(2)$^{\cite{BM}}$ \ {\it $\mathcal{D}^{+}(A)$ is equivalent, as a
right triangulated category, to the factor module category
$A^\backprime $-${\rm \overline{mod}}$.}

\vskip 0.2in

{\bf Lemma 2.2.}$^{\cite{HW}}$\ {\it Let $M$ be an indecomposable
$\hat{A}$-module which is not projective-injective. Then there
exists an indecomposable $A$-module $N$ such that $M\simeq
\Omega^{-l}_{\hat{A}}N$ for some $l\in \mathbb{Z}$. We denote by $l$
the degree of $M$, that is, ${\rm deg} M=l$. By abuse language, we
also call ${\rm deg} M$ the degree of $M$ in $\mathcal{D}^{b}(A)$.}

\vskip 0.2in

Let $\mathscr{A}$ be a triangulated category with shift functor [1].
Let $\mathcal {X}, \mathcal {Y}$ be two subcategories of
$\mathscr{A}$. Then we denote by $\mathcal {X}\ast\mathcal {Y}$ the
subcategory of $\mathscr{A}$ consisting of $M$ with a triangle
$X\rightarrow M\rightarrow Y\rightarrow X[1]$ such that
$X\in\mathcal {X}$ and $Y\in \mathcal {Y}$.

\vskip 0.2in

Throughout this paper, we follow the standard terminology and
notation used in the representation theory of algebras, see
\cite{ASS, ARS} and \cite{H, R}.

\vskip 0.2in

\section {Endomorphism algebras arising from silting mutation}

\vskip 0.2in

Let $A$ be a hereditary finite-dimensional algebra with $n$
isomorphism classes of simple modules, and $\mathcal {D} = D^b(A)$
be the bounded derived category of finitely generated $A$ modules.
Recall that for an integer $m \geq 1$, the $m$-cluster category of
$A$ is the orbit category $\mathscr{C}_m = \mathcal
{D}/\tau^{-1}[m]$, where $\tau$ is the AR-translation in $\mathcal
{D}$, see \cite{HT, BZ}, and $[m]$ is the $m$-fold composition of
$[1]$, the shift functor. This category is known to be triangulated
by \cite{K}, and  also a Krull-Schmidt category.

\vskip 0.2in

A basic object $T$ in $\mathcal {D}$ is said to be partial silting
if ${\rm Ext}^i(T , T ) = 0$ for $i > 0$, and silting if in addition
$T$ is maximal with this property.

\vskip 0.2in

Let $T$ be a basic silting object in $\mathcal {D}$. Without loss of
generality we can assume that $T$ is in the non-negative part
$\mathcal {D}^+$  of $\mathcal {D}$, i.e. we assume
$$
T = T_0[0] \oplus T_1[1]\oplus \cdots \oplus T_m[m] \ \ \ \ \ \ \ \
\ \ \ \ \ \ \ \ (1)$$ with each $T_i$ in ${\rm mod}~A$.

We fix
$$
\mathcal {S}_m = {\rm mod}~A[0] \vee {\rm mod}~A[1]\vee \cdots \vee
{\rm mod}~A[m-1] \vee A[m].
$$
Then $\mathcal {S}_m$ is a fundamental domain for $\mathscr{C}_m$ in
$\mathcal {D}$; this means that the map from isomorphism classes of
objects in $\mathcal {S}_m$ to isomorphism classes of objects in
$\mathscr{C}_m$ is bijective.

\vskip 0.2in

{\bf Lemma 3.1.}$^{\cite{BRT}}$ \  {\it Let $A$ be a
finite-dimensional hereditary algebra, and let $\mathcal {S}_m$ be
the fundamental domain as above for the m-cluster category
$\mathscr{C}_m$ of $A$. Let $T$ be an object in $\mathcal {S}_m$.
Then we have the following:

\vskip 0.1in

(a) $T$ is a partial silting object in $\mathcal {D}$ if and only if
$T$ is rigid in $\mathscr{C}_m$.

\vskip 0.1in

(b) $T$ is a silting object in $\mathcal {D}$ if and only if $T$ is
an m-cluster tilting object in $\mathscr{C}_m$.}

\vskip 0.2in

{\bf Theorem 3.2.}\  {\it Let $T$ be a basic silting object in
$\mathcal{S}_m$, and $\Gamma= {\rm End}_{\mathcal{D}}T$. Let
$\mathscr{T}= {\rm add}~T$ and $\mathcal {W}=\mathscr{T}\ast
\mathscr{T}[1]$. Then the functor
$$G = {\rm Hom}_{\mathcal {D}}(T, -
): \mathcal {W} \rightarrow {\rm mod}~\Gamma $$ is full and dense,
and $G$ also induces an equivalence functor
$$\overline{G} : \mathcal {W} /{\rm add}~(T[1]) \rightarrow  {\rm
mod}~\Gamma.$$}

\vskip 0.1in

{\bf Proof.}\ \ First, $G$ is dense.

Indeed, take $X\in{\rm add}~\Gamma$, let
$P_1\stackrel{\alpha}\longrightarrow P_0\rightarrow X\rightarrow 0$
be a minimal projective resolution of $X$. Then there exists a
morphism $f: T_1\rightarrow T_0$ with $T_1, T_0\in {\rm add}~T$ such
that $G(f)= \alpha$, $G(T_1)=P_1$ and $G(T_0)=P_0$. Thus we obtain a
triangle in $\mathcal {D}$
$$T_1\stackrel{f}\longrightarrow T_0\rightarrow Y\rightarrow
T_1[1],$$ applying $G = {\rm Hom}_{\mathcal {D}}(T, - )$ yields an
exact sequence
$$
G(T_1)\stackrel{G(f)}\longrightarrow G(T_0)\rightarrow
G(Y)\rightarrow G(T_1[1])=0,
$$
which implies that $G(Y)\simeq X$, that is  $G$ is dense.

\vskip 0.1in

Now, we show that $G$ is full. For any $h\in{\rm Hom}_\Gamma (X,
Y)$, there exist $M$ and $N$ in $\mathcal {W}$ with $X= G(M)$ and
$Y=G(N)$, since $G$ is dense.

We have triangles $T_1\stackrel{f}\longrightarrow
T_0\stackrel{\pi}\longrightarrow M\rightarrow T_1[1]$ and
$T_1'\stackrel{f'}\longrightarrow T_0'\stackrel{\pi'}\longrightarrow
N\rightarrow T_1[1]$ with $T_i, T_i'\in {\rm add}~T$ for $i=0, 1$.
Applying $G$ yields a commutative diagram with exact rows in ${\rm
mod}~\Gamma$
$$
\begin{array}{rclcl}
 G(T_1)&\stackrel{G(f)}{\longrightarrow}&  G(T_0)
 &\stackrel{G(\pi)}{\longrightarrow} & X \rightarrow 0 \\
 h_1 \downarrow & &\downarrow h_0&& \downarrow h\\
  G(T_1')& \stackrel{G(f')}{\longrightarrow} &  G(T_0')
&\stackrel{G(\pi')}{\longrightarrow}& Y \rightarrow 0 .
\end{array}
$$
Note that $G: {\rm add}~T \rightarrow {\rm add}~_\Gamma\Gamma$ is an
equivalence functor, and that $G(T_i)$ and $G(T_i')$ are projective
$\Gamma$ modules for $i=0, 1$. Thus we get the following commutative
diagram:
$$
\begin{array}{rcl}
T_1&\stackrel{f}{\longrightarrow}& T_0 \\
 g_1 \downarrow & &\downarrow g_0\\
 T_1'& \stackrel{f'}{\longrightarrow} &  T_0' .
\end{array}
$$
Then we have the following commutative diagram of triangles
$$\begin{array}{rclcll}
 T_1 &\stackrel{f}{\longrightarrow}&  T_0
 &\stackrel{ \pi }{\longrightarrow} & M \rightarrow &T_1[1] \\
 g_1 \downarrow & &\downarrow g_0& & \downarrow g  &\downarrow g[1] \\
  T_1' & \stackrel{ f' }{\longrightarrow} & T_0'
&\stackrel{\pi'}{\longrightarrow}& N \rightarrow &T_1'[1] ,
\end{array}
$$
thus $G(g)=h$, namely $G$ is full.

\vskip 0.1in

Finally, let $\lambda: M\rightarrow N$ in $\mathcal {W}$ with
$G(\lambda)=0$, that is, ${\rm Hom}_{\mathcal {D}}(T,\lambda): {\rm
Hom}_{\mathcal {D}}(T, M) \rightarrow  {\rm Hom}_{\mathcal {D}}(T,
N)$ is 0. Consider the diagram
$$\begin{array}{rclcll}
 T_1 &\stackrel{f}{\longrightarrow}&  T_0
 &\stackrel{ \pi }{\longrightarrow} & M \stackrel{\sigma}{\longrightarrow}  &T_1[1] \\
  & && & \downarrow \lambda  & \\
  &  & & & N  &
\end{array}
$$
where $T_1 \stackrel{f}{\longrightarrow}  T_0  \stackrel{ \pi
}{\longrightarrow}  M \stackrel{\sigma}{\longrightarrow}  T_1[1]$ is
a triangle.
 Since $\lambda \pi = 0$, there is a map $t: T_1[1] \rightarrow  N$
such that $t\sigma=\lambda$, that is,  $\lambda: M\rightarrow N$
factors through ${\rm add}~(T[1])$ . This shows that $\overline{G}$
is faithful, and consequently an equivalence.           $\hfill\Box$

\vskip 0.2in

{\bf Remark.}\ This theorem implies that we can obtain many abeliean
categories from a triangulated category $\mathcal {D}$, which is of
independent interests.

\vskip 0.2in

Let $\mathscr{E}$ be the set of all basic silting objects in
$\mathcal {D}^+$. We define the silting quiver
$\overrightarrow{\mathscr{S}}_{\mathscr{E}}$ of $\mathcal {D}^+$ as
follows.

The vertices of $\overrightarrow{\mathscr{S}}_{\mathscr{E}}$ is the
elements in $\mathscr{E}$, and for two elements $W_1, W_2$ in
$\mathscr{E}$, there is an arrow $W_1\rightarrow W_2$  if and only
if there exists an almost silting object $T$ with $W_1=T\oplus X$
and $W_2=T\oplus Y$ such that there exists a triangle
$$
(*)\ \ \ \ \ \ \ \ \ \ \ \  X \stackrel{f}\longrightarrow
B\stackrel{g}\longrightarrow Y\rightarrow X[1]
$$
with $f$ ($g$) being a left (right) minimal ${\rm
add}~T$-approximation. $(*)$ is called mutation triangle.

\vskip 0.2in

{\bf Proposition 3.3.}\ {\it Silting quiver
$\overrightarrow{\mathscr{S}}_{\mathscr{E}}$ is connected.}

\vskip 0.1in

{\bf Proof.}\ \  Let $W_1$ and $W_2$ be two elements of
$\overrightarrow{\mathscr{S}}_{\mathscr{E}}$. Without loss of
generality we can assume that $W_1$ and $W_2$ lie in $\mathcal
{S}_m$ for some positive integer $m$. According to Theorem 3.5 in
\cite{BRT}, $W_1$ and $W_2$ are tilting objects in $m$-cluster
category $\mathscr{C}_m$.  By Proposition 4.5 in \cite{ZZ} we know
that the $m$-cluster tilting quiver of $\mathscr{C}_m$ is connected,
and by using Theorem 3.5 in \cite{BRT} again, we know  that
$\overrightarrow{\mathscr{S}}_{\mathscr{E}}$ is connected.  This
completes the proof. $\hfill\Box$

\vskip 0.3in

Let $T = T_1\oplus \cdots\oplus T_n$ be a basic silting object in
$\mathcal {D} = D^b(A)$, where the $T_i$ are indecomposable and $n$
is the number of isomorphism classes of simple $A$-modules. We
assume without loss of generality that $T$ is in $\mathcal {D}^+$.
Let $m$ be an integer such that $T$ is in the fundamental domain
$\mathcal {S}_m$ of the $m$-cluster category $\mathscr{C}_m$. Fix an
indecomposable direct summand $T_i$, with $1 \leq i \leq n$.  Recall
from \cite{BRT}, we know that the almost complete silting object
$T/T_i$ has a countably infinite number of non-isomorphic
complements $M_i$ for $i\in \mathbb{Z}$. In particular, there are
complements $M_{-1}$ and $M_{m+1}$, such that $M_j\simeq M_{-1}[j +
1]$ for $j < -1$, and $M_j = M_{m+1}[j-(m +1)]$ for $j
>m +1$.

\vskip 0.2in

{\bf Theorem 3.4.}\ {\it Let $T = T_1\oplus \cdots\oplus T_n$ be a
basic silting object in $\mathcal {S}_m$. Fix an indecomposable
direct summand $T_i$, with $1 \leq i \leq n$.  For each silting
complements $M_j$ of $T/T_i$, we set $\Gamma_j= {\rm End}_{\mathcal
{D}}(M_j\oplus T/T_i)$ and  $\mathcal {W}_j=\mathscr{T}_j\ast
\mathscr{T}_j[1]$  with $\mathscr{T}_j={\rm add}~(M_j\oplus T/T_i)$.
Then

\vskip 0.1in

(1)\ The $\Gamma_j$ module ${\rm Hom}_{\mathcal {D}}(M_j\oplus
T/T_i, M_{j-1}[1])$ is simple.

\vskip 0.1in

(2)\  Let $S_{M_j}$ be the simple top of the indecomposable
projective $\Gamma_j$ module ${\rm Hom}_{\mathcal {D}}(M_j\oplus
T/T_i, M_{j})$. Then we have
$$
{\rm mod}~\Gamma_j /{\rm add}~S_{M_j} \simeq \mathcal {W}_j/{\rm
add}~( M_{j}\oplus M_{j-1}\oplus T/T_i)[1].
$$
Moreover, if $j>m$, then $\Gamma_j \simeq \left (
\begin{array}{cc}
{\rm End}_{\mathcal {D}}M_j & 0\\
0 & {\rm End}_{\mathcal {D}}(T/T_i) \end{array}\right )$.}

\vskip 0.1in

{\bf Proof.}\   (1) \  ${\rm Hom}_{\mathcal {D}}(M_j,
M_{j-1}[1])\simeq {\rm Ext}^1_{\mathcal {D}}(M_j, M_{j-1})$ is
one-dimensional over the factor algebra ${\rm End}_{\mathcal
{D}}(M_j)/ {\rm Rad}(M_j,M_j)$, see \cite{ABST, KR}, and thus a
simple ${\rm End}_{\mathcal {D}}~(M_j)$ module. We have
\begin{eqnarray*}
{\rm Hom}_{\mathcal {D}}(M_j\oplus T/T_i, M_{j-1}[1])&\simeq & {\rm
Hom}_{\mathcal {D}}(T/T_i, M_{j-1}[1])\oplus {\rm Hom}_{\mathcal
{D}}(M_j, M_{j-1}[1])\\
&\simeq & {\rm Hom}_{\mathcal {D}}(M_j, M_{j-1}[1])
\end{eqnarray*}
Thus,  ${\rm Hom}_{\mathcal {D}}(M_j\oplus T/T_i, M_{j-1}[1])\simeq
{\rm Hom}_{\mathcal {D}}(M_j, M_{j-1}[1])$ as a  $\Gamma_j$ module,
and is hence simple.

\vskip 0.1in

(2)\  Let $M_{j-1} \longrightarrow B\stackrel{f}\longrightarrow
M_j\rightarrow M_{j-1}[1]$  be a triangle with $f$ being a minimal
right ${\rm add}~(T/T_i)$-approximation. Then $M_{j-1}$ is a silting
complement of $T/T_i$, and $M_{j-1}[1]\in \mathcal {W}_j$. Applying
${\rm Hom}_{\mathcal {D}}(M_j\oplus T/T_i, -)$ yields an exact
sequence
$$
{\rm Hom}_{\mathcal {D}}(M_j\oplus T/T_i, B)\rightarrow {\rm
Hom}_{\mathcal {D}}(M_j\oplus T/T_i, M_j)\rightarrow {\rm
Hom}_{\mathcal {D}}(M_j\oplus T/T_i, M_{j-1}[1]) \rightarrow 0,
$$
which implies that $S_{M_j}\simeq {\rm Hom}_{\mathcal {D}}(M_j\oplus
T/T_i, M_{j-1}[1])$.

According to Theorem 3.2, we get an equivalence $\mathcal {W}_j/
{\rm add}~(M_j\oplus T/T_i)[1]\rightarrow {\rm mod}~ \Gamma_j$ such
that $ M_{j-1}[1]\mapsto S_{M_j}$. Hence we get the following.
$$
{\rm mod}~\Gamma_j /{\rm add}~S_{M_j} \simeq \mathcal {W}_j/{\rm
add}~( M_{j}\oplus M_{j-1}\oplus T/T_i)[1].
$$

Note that
$$
\Gamma_j={\rm End}_{\mathcal {D}}(M_j \oplus T/T_i)=  \left (
\begin{array}{cc}
{\rm End}_{\mathcal {D}}M_j & {\rm Hom}_{\mathcal {D}}(M_j,
T/T_i) \\
{\rm Hom}_{\mathcal {D}}(T/T_i, M_j) & {\rm End}_{\mathcal
{D}}(T/T_i) \end{array}\right )
$$
and that $T/T_i\in \mathcal {S}_m$, thus ${\rm Hom}_{\mathcal
{D}}(M_j, T/T_i)=0$ for all $j> m$, since ${\rm deg}~M_j> m$.

Assume  $j=m+1$ and let $T/T_i=X_m\oplus P[m]$ with $P$ projective
and ${\rm deg}~X_m\leq  m-1$. Then
$$
{\rm Hom}_{\mathcal {D}}(T/T_i, M_{m+1})={\rm Hom}_{\mathcal
{D}}(X_m, M_{m+1})\oplus {\rm Hom}_{\mathcal {D}}(P[m], M_{m+1})=0.
$$

If $j>m+1$, then ${\rm Hom}_{\mathcal {D}}(T/T_i, M_{j})=0$, since
$T/T_i\in \mathcal {S}_m$ and ${\rm deg}~M_j\geq m+2$.

Hence for all $j> m$, we have $\Gamma_j \simeq \left (
\begin{array}{cc}
{\rm End}_{\mathcal {D}}M_j & 0\\
0 & {\rm End}_{\mathcal {D}}(T/T_i) \end{array}\right )$.  The proof
is completed.       $\hfill\Box$

\vskip 0.2in

{\bf Remark.}\  We do not know whether these endomorphism algebras
$\Gamma_j$ for $j\geq 0$ are derived equivalent.

\vskip 0.3in

 Let $M = M_1\oplus M_2\oplus\cdots\oplus M_n$, where the $M_i$ are
pairwise non-isomorphic indecomposable objects in $\mathcal {D}$.

Let $S_i = S_{M_i}$ be the simple ${\rm End}_{\mathcal
{D}}(M)$-module corresponding to $M_i$, and let $P_i = {\rm
Hom}_{\mathcal {D}} (M, M_i)$ be the indecomposable projective ${\rm
End}_{\mathcal {D}}M$ module with top $S_i$.

\vskip 0.2in

Let $Q_M$ be the quiver of $\Gamma_M={\rm End}_{\mathcal {D}}(M)$.
For $1\leq i,j \leq n$, the following numbers are equal:

\vskip 0.1in

$\bullet$ \  The number of arrows $i\rightarrow j$ in the quiver
$Q_M$ (Note that ${\rm End}_{\mathcal {D}} M\simeq k\ Q_M/I$);

\vskip 0.1in

$\bullet$ \ ${\rm dim}~{\rm Ext}^1_{\Gamma_M}(S_i,S_j)$;

\vskip 0.1in

$\bullet$ \ The dimension of the space of irreducible maps
$M_i\rightarrow M_j$ in the category ${\rm add} (M)$;

\vskip 0.1in

$\bullet$ \ The dimension of the space of irreducible maps
$P_i\rightarrow P_j$ in the category ${\rm add}\
(P_1\oplus\cdots\oplus P_n)$ of projective ${\rm End}_{\mathcal
{D}}M$ modules.

\vskip 0.1in

Furthermore, let $f: M_i\rightarrow M'$ (resp. $g: M'' \rightarrow
M_i$ ) be a minimal left (resp. right) ${\rm add}\
(M/{M_i})$-approximation of $M_i$. If $i\neq j$, then we have

$$ {\rm dim}~{\rm Ext}^1_{\Gamma_M}(S_i,S_j)= [M': M_j],$$

$$ {\rm dim}~{\rm Ext}^1_{\Gamma_M}(S_j,S_i)= [M'' : M_j].$$

\vskip 0.2in

{\bf Proposition 3.5.}\ {\it Let $T=T_1\oplus T_2\oplus\cdots\oplus
T_n$ be a basic silting object in $\mathcal{S}_m$, and let  $Q_T$
 be the quiver of $\Gamma_T= {\rm End}_{\mathcal{D}}~T$. Then

\vskip 0.1in

(1) \ $Q_T$  has no loops and no 2-cycles.

\vskip 0.1in

(2) \  $\Gamma_T$ is quasi-hereditary. In particular, ${\rm
gl.dim}~\Gamma_T < \infty$.}

\vskip 0.1in

{\bf Proof.}\ (1)\ Since $T$ can be regarded as a cluster tilting
object in $\mathscr{C}_m$, and $Q_T$ is a subquiver of the quiver
$Q_{{\rm End}_{\mathscr{C}_m}T}$ of the endomorphism algebra ${\rm
End}_{\mathscr{C}_m}T$ with same vertices, which implies that $Q_T$
has no loops and no 2-cycles.

\vskip 0.1in

(2) We can write $T$ as $T=M_0\oplus M_1[1]\oplus \cdots\oplus
M_{m-1}[m-1]\oplus P[m]$, where $M_1,\cdots, M_{m-1}$ and $P$ are
$A$ modules with $P$ projective.   Then $M_0,\cdots, M_{m-1}, P$ are
partial tilting $A$ modules, which implies $Q_T$ has no cycles,
hence $\Gamma_T$ is quasi-hereditary, see \cite{R1}. Hence, ${\rm
gl.dim}~\Gamma_T < \infty$.
   $\hfill\Box$

\vskip 0.2in

{\bf Remark.}\  The situation is very different in $m$-cluster
category $\mathscr{C}_m$. For example, the quiver $Q_ {{\rm
End}_{\mathscr{C}_m}~T}$ of ${\rm End}_{\mathscr{C}_m}~T$ usually
has cycles for a cluster tilting object in $\mathscr{C}_m$.

\vskip 0.2in

\section {Endomorphism algebras determined by tilting mutations in $A^{(m)}$}

Let $\Lambda$ be a finite dimensional algebra over field $k$, and
let $S$ be a non-injective simple $\Lambda$ module with the
following two properties:

(a) ${\rm proj.dim}_{\Lambda}(\tau^{-1}\ S) \leq 1$, and

\vskip 0.07in

(b) ${\rm Ext}^1_{\Lambda}(S, S) = 0$.

\vskip 0.2in

Here $\tau^{-1}= {\rm Tr D}$ stands for the inverse of the
Auslander-Reiten translation.  We denote the projective cover of $S$
by $P(S)$, and assume that $\Lambda = P(S) \oplus P$ such that there
is not any direct summand of $P$ isomorphic to $P(S)$. Let $T =
\tau^{-1} S \oplus P$. It is well known that $T$ is a tilting
module. Such a tilting module is called a BB-tilting module.
Moreover, if $S$ is also a projective non-injective simple module,
then ${\rm Hom}_{\Lambda}(D(\Lambda),S) = 0$, and therefore ${\rm
proj.dim}_{\Lambda}(\tau^{-1}\ S) \leq 1$. Thus $T$ is a BB-tilting
module. This special tilting module  is called an APR-tilting module
in literature. It is widely used in the representation theory of
algebras. Note that if $S$ is a non-injective, projective simple
$\Lambda$-module, then there is an Auslander-Reiten sequence
$0\rightarrow S\rightarrow P' \rightarrow \tau^{-1}\ S \rightarrow
0$ in $\Lambda$-mod with $P'$ projective.

\vskip 0.2in

Let $M$ be a faithful almost tilting $\Lambda$-module, and let $X$
be an indecomposable complement to $X$ which is cogenerated by $M$.
According to \cite{HU1, HU2}, we know that there exists an exact
sequence
$$
(\dag)\ \ \ \ \ \ \ \ \ \ \ \ \ \ 0\rightarrow X\stackrel{f}
\longrightarrow M'\stackrel{g} \longrightarrow Y \rightarrow 0
$$
with $Y$ indecomposable and $f$
($g$) being a minimal left (right) ${\rm add}\ M$-approximation.

\vskip 0.2in

{\bf Lemma 4.1.}\ {\it Take the notations as above. Let $\Gamma={\rm
End}_\Lambda\ (X\oplus M)$. Then there exists a BB-tilting $\Gamma$
module $T= \tau^{-1} S\oplus P$ such that ${\rm End}_{\Gamma}\
T\simeq \ {\rm End}_{\Lambda}\ (M\oplus Y)$.}

\vskip 0.1in

{\bf Proof}\ \  Let $V=X\oplus M$. Applying ${\rm
Hom}_{\Lambda}(V,-)$ to $(\dag)$ yields an exact sequence
$$
(\ddag)\ \ \ \ \ \ \  0\rightarrow {\rm
Hom}_{\Lambda}(V,X)\stackrel{f^*} \longrightarrow {\rm
Hom}_{\Lambda}(V,M') \longrightarrow L \rightarrow 0 .
$$
Note that $(\ddag)$ is the minimal projective resolution of $L$. Let
$T=L\oplus {\rm Hom}_{\Lambda}(V,M)$. Then $T$ is a tilting module
by the proof of Lemma 3.4 in \cite{HX}.

Applying ${\rm Hom}_{\Gamma}(-,\Gamma)$ to $(\ddag)$ we get an exact
sequence of $\Gamma$ modules
$$
{\rm Hom}_{\Gamma}({\rm Hom}_{\Lambda}(V,M'),\Gamma) \longrightarrow
{\rm Hom}_{\Gamma}({\rm Hom}_{\Lambda}(V,X),\Gamma) \longrightarrow
{\rm Tr}_{\Gamma}\ L \rightarrow 0,
$$
which is isomorphic to the following exact sequence
$$
{\rm Hom}_{\Lambda}(M',V)\stackrel{f_*} \longrightarrow {\rm
Hom}_{\Lambda}(X,V) \longrightarrow {\rm Tr}_{\Gamma}\ L\rightarrow
0,
$$
where $f_*={\rm Hom}_{\Lambda}(f, V)$.

We claim that ${\rm Im}\ f_*$ is the radical of the indecomposable
projective $\Gamma$ module ${\rm Hom}_{\Lambda}(X,V)$.

Indeed,  $f_*={\rm Hom}_{\Lambda}(f, X\oplus M)={\rm
Hom}_{\Lambda}(f, X)\oplus {\rm Hom}_{\Lambda}(f, M)$.

Since $f$ is a minimal left ${\rm add}\ M$ approximation of $X$, by
applying ${\rm Hom}_{\Lambda}(-,M)$ to $(\dag)$ we get an exact
sequence
$$
0\rightarrow {\rm Hom}_{\Lambda}(Y,M) \longrightarrow {\rm
Hom}_{\Lambda}(M',M) \stackrel{{\rm Hom}_{\Lambda}(f,
M)}\longrightarrow {\rm Hom}_{\Lambda}(X,M)\rightarrow 0.
$$
Hence ${\rm Hom}_{\Lambda}(f, M)$ is surjective.

By applying ${\rm Hom}_{\Lambda}(-,X)$ to $(\dag)$ we have an exact
sequence
$$
{\rm Hom}_{\Lambda}(M',X) \stackrel{{\rm Hom}_{\Lambda}(f,
X)}\longrightarrow {\rm Hom}_{\Lambda}(X,X) \longrightarrow {\rm
Ext}_{\Lambda}^1(Y,X)\rightarrow 0,
$$
it forces that ${\rm Im~Hom}_{\Lambda}(f,X)= {\rm rad}\ {\rm
Hom}_{\Lambda}(X,X)$ since ${\rm dim}_k\ {\rm Ext}_{A}^1(Y,X)=1$. It
follows that ${\rm Im}\ f_*={\rm rad}\ {\rm Hom}_{A}(X,V)$, and our
claim is true.

It follows from our claim that ${\rm Tr}_{\Gamma}\ L$ is a simple
$\Gamma$ module and $\tau L_{\Gamma}$ is the simple socle $S$ of the
indecomposable injective $\Gamma$ module $D{\rm
Hom}_{\Lambda}(X,V)$, and we know that $L\simeq \tau_{\Gamma}^{-1}
S$. Then $T$ is a BB-tilting $\Gamma$ module. By Lemma 3.4 in
\cite{HX} again, we have that ${\rm End}_{\Gamma}\ T \simeq {\rm
End}_{\Lambda}\ (M\oplus Y)$. $\hfill\Box$

\vskip 0.2in

Now, let $A$ be a hereditary algebra over an algebraically closed
field $k$. The $m$-replicated algebra $A^{(m)}$ of $A$ is defined as
the quotient of  the repetitive algebra $\hat{A}$, that is,
$$A^{(m)}=\begin{pmatrix}
                  A_{0} &  &  &  & 0 & \\
                  Q_{1} & A_{1} & &  &  &  \\
                   & Q_{2} & A_{2} &  &  &  \\
                 &  &  & \ddots & \ddots &  \\
                  &  0&  &  & Q_{m} & A_{m} \\
                \end{pmatrix}.
$$

\vskip 0.2in

Let $M$ be a faithful almost complete tilting $A^{(m)}$-module.
Recall from \cite{LZ, Z} we know that $M$ has non-isomorphic
indecomposable complements $X_0, \cdots, X_t$ with  $2m\leq t\leq
2m+1$.

\vskip 0.2in

{\bf Theorem 4.2.} \ {\it Let $M$ be a faithful almost complete
tilting $A^{(m)}$-module with non-isomorphic indecomposable
complements $X_0, \cdots, X_t$, and let $\Gamma_i={\rm
End}_{A^{(m)}}(X_i\oplus M)$. Then there is a BB-tilting $\Gamma_i$
module $T_i$ such that ${\rm End}_{\Gamma_i} T_i\simeq \Gamma_{i+1}$
for $0\leq i \leq t-1$.}

\vskip 0.1in

{\bf Proof}\ \ According to \cite{LLZ, LZ}, we have an exact
sequence
$$
0\longrightarrow X_i \stackrel {f_i}\longrightarrow B_i \stackrel
{g_i}\longrightarrow X_{i+1}\longrightarrow 0
$$
such that  $f_i$ is a minimal left ${\rm add}~M$-approximation of
$X_i$ and that $g_i$ is a minimal right ${\rm add}~M$-approximation
of $X_{i+1}$. By Lemma 4.1, there exists a BB-tilting  $\Gamma_i$
module $T_i$ such that ${\rm End}_{\Gamma_i} T_i\simeq
\Gamma_{i+1}$.              $\hfill\Box$

\vskip 0.2in

{\bf Theorem 4.3.}\  {\it Let $T= T_1\oplus T_2\oplus\cdots\oplus
T_n$ be a basic silting object in $\mathcal{S}_m$, and let  $P$ be
the direct sum of all indecomposable injective-projective $A^{(m)}$
modules. Then $\Gamma_T= {\rm End}_{\mathcal{D}}~T\simeq {\rm
\overline{End}}_{A^{(m)}}~(T\oplus P)$.}

\vskip 0.1in

{\bf Proof.}\ ${\rm \overline{End}}_{A^{(m)}}~(T\oplus P)={\rm
End}_{A^{(m)}}~(T\oplus P)/\mathscr{P}\simeq {\rm
End}_{\mathcal{D}^+}~T\simeq {\rm End}_{\mathcal{D}}~T$.
$\hfill\Box$

\end{document}